\documentclass[a4paper]{article}
\usepackage{epsfig,amsmath,amsfonts,amsthm,mathrsfs}
\usepackage{verbatim}

\title{How to find horizon-independent optimal strategies leading off to
infinity: a max-plus approach}

\author{Marianne Akian and St\'ephane Gaubert and Cormac Walsh
\thanks{The authors are with INRIA, Domaine de Voluceau,
78153 Le Chesnay C\'edex, France. 
{\tt Marianne.Akian@inria.fr},
{\tt Stephane.Gaubert@inria.fr},
{\tt Cormac.Walsh@inria.fr}}%
}

\usepackage{hyperref}
\newcommand\newgeo{\gamma'}

\newcommand\geo{\gamma}
\newcommand{\states}{X}
\newcommand{\recclasses}{\overline{\states}}
\newcommand{\rplus}{\R_+}
\newcommand{\sS}{\mathscr{S}}
\newcommand{\sH}{\mathscr{H}}
\newcommand{\R}{\mathbb{R}}
\newcommand{\N}{\mathbb{N}}
\newcommand{\sB}{\mathscr{B}}

\newcommand{\sK}{\mathscr{K}}
\newcommand{\sM}{\mathscr{M}}
\newcommand{\sMin}{\sM^{m}}
\newcommand{\set}[2]{\{#1\mid #2\}}
\def\of#1,#2{(#1,#2)}

\newcommand\minspace{\sMin}
\newcounter{change}

\def\bfx{\mathbf{x}}
\def\bfu{\mathbf{u}}
\def\new#1{{\em #1}}
\newcommand{\mrm}[1]{\text{\rm #1}}
\newcommand{\rmax}{\R_{\max}}
\newcommand{\rmaxb}{\overline{\R}_{\max}}

\newtheorem{prop}{Proposition}[section]
\newtheorem{proposition}[prop]{Proposition}
\newtheorem{corollary}[prop]{Corollary}
\newtheorem{lemma}[prop]{Lemma}
\newtheorem{assumption}[prop]{Assumption}

\newtheorem{theorem}[prop]{Theorem}
\theoremstyle{definition}
\newtheorem{definition}[prop]{Definition}
\theoremstyle{remark}
\newtheorem{example}[prop]{Example}
\newtheorem{remark}[prop]{Remark}
\begin{document}

\maketitle

{To appear in Proceedings of the 45th IEEE Conference on Decision and Control, 2006.} \\

\begin{abstract}
A general problem in optimal control consists of finding
a terminal reward that makes the value function independent
of the horizon. Such a terminal reward can be interpreted
as a max-plus eigenvector of the associated Lax-Oleinik semigroup.
We give a representation formula for all these eigenvectors, which
applies to optimal control problems in which the state space is non compact.
This representation involves an abstract boundary of the state space,
which extends the boundary of metric spaces defined
in terms of Busemann functions (the horoboundary).
Extremal generators of the eigenspace correspond to certain
boundary points, which are the limit of almost-geodesics.
We illustrate our results in the case of a linear quadratic problem.
\end{abstract}

\section{Introduction}
A basic optimal control problem with finite horizon consists of maximising
a reward of the form
\begin{subequations}\label{pb1}
\begin{align}\label{e-reward}
\int_0^t L(\bfx(s),\bfu(s)) d s + \phi(\bfx(t)) \enspace,
\end{align}
over the trajectories $\bfx(\cdot)$ and $\bfu(\cdot)$ satisfying
\begin{align}
\dot{\bfx}(s) &= f(\bfx(s),\bfu(s)) \enspace,  \label{e-dyn}\\
\bfu(s)& \in U,\quad
\bfx(s) \in \states
\qquad \text{for all $0\leq s\leq t$}, \label{e-set}\\
\text{and}\qquad
 \bfx(0) &= x \enspace .\label{e-init}
\end{align}\end{subequations}
We are given the state space $\states\subset \R^n$,
the set of admissible controls $U\subset \R^p$,
the dynamics $f: \states\times U \to \R^n$,
the initial condition $x\in \R^n$,
the horizon $t$, the Lagrangian $L: \states\times U \to \R$,
and the terminal reward $\phi: \states\to \R\cup\{-\infty\}$.
Of course, the data of the problem must satisfy certain
technical assumptions that we skip for the moment,
and the trajectories $\bfx(\cdot)$ and $\bfu(\cdot)$ must 
have a certain regularity: $\bfu$ is required to be measurable
and $\bfx$ is required to be absolutely continuous. We are interested in
the maximal value of the reward, which we denote by $v^t(x)$, 
and in the optimal solutions $\bfx(\cdot)$ and $\bfu(\cdot)$.

One may wish in particular to determine
a final reward $\phi$ making the \new{value function} 
$v^t: \states\to \R\cup\{-\infty\}$, $x\mapsto v^t(x)$,
independent of the horizon $t$ in a reasonable 
sense. Since $v^t$ is expected to grow linearly with the horizon $t$,
we look for final rewards $\phi$ such that
\begin{align}
v^t = \lambda t + \phi \qquad \text{for all $t \geq 0$}\enspace ,
\label{ergo}
\end{align}
where $\lambda$ is a constant that can be interpreted
as a mean reward per time unit. (For any function
$\phi$ and constant $\mu$, we denote by $\mu+\phi$ the function
sending $x$ to $\mu+\phi(x)$.)

In this paper, we address the problem of finding all the solutions $\phi$
of~\eqref{ergo} corresponding to a given constant $\lambda$.

The economic motivation of this classical problem is
presented in~\cite{yakovenko}:
in finite horizon optimal control, the optimal strategy may change
dramatically as the current time gets closer to the horizon,
which is often undesirable.
To give a simple illustration, 
imagine a company having a concession
to run a network belonging to a city during a given period.
The company might stop repairing the network
as time approaches the horizon.
Thus, the city might wish to choose a terminal penalty $\phi$ 
preventing the company's strategy from being distorted by the horizon.
A way of addressing this issue is to solve~\eqref{ergo}.

To study this problem, it is convenient
to introduce the \new{Lax-Oleinik semigroup}, $(S^t)_{t\geq 0}$, 
which consists of the operators $S^t$ sending the terminal reward $\phi$
to the value function $v^t$. 
Observe that the operator $S^t$ is \new{max-plus linear},
meaning that it preserves
suprema and commutes with the addition of a constant.
Equation~\eqref{ergo} can be rewritten as
\begin{align}
S^t\phi = \lambda t + \phi \qquad \text{for all $t \geq 0$} \enspace.
\label{spectral}
\end{align}
We recognise this as a max-plus analogue of the spectral equation.
We say that $\phi$ is an \new{eigenfunction}
of the semigroup and that $\lambda$ is the associated \new{eigenvalue}.

The discrete-time analogue of problem~\eqref{spectral} has received
a considerable amount of attention~\cite{bcoq,bapat98,gondran02,AGW-s,AGW-m}. 
Problem~\eqref{spectral} itself has been
studied by Maslov, Kolokoltsov, and other members
of the ``idempotent analysis'' school~\cite{kolokoltsov_maslov,maslov92}. 
More recently, it has appeared in the setting 
of Fathi's ``weak KAM'' theory~\cite{fathi97b,fathi97c,fathi03}.
Results of idempotent analysis and of weak KAM theory yield
a representation of all the solutions of~\eqref{spectral}, 
under certain assumptions. For instance, the approach
of~\cite{fathi03} is in the setting where $\states$ is a Riemannian manifold
and the Lagrangian has smoothness and strict-convexity
properties. The representation results
of~\cite[Section~3.2]{kolokoltsov_maslov}, on the other hand, apply to other
special cases.
However, the discrete-time theory~\cite{AGW-m} suggests that 
general results should hold in the continuous-time setting.
Such a generality is required in many applications. For instance,
the Lagrangian might not be smooth
or the trajectory $\mathbf{x}$ may be the sum of an absolutely
continuous part and a singular part modelling ``jumps''.

In~\cite{AGW-m}, some continuous-time results were obtained in special cases
by applying the discrete-time results.
However, we obtain more general results when we redo the theory from the
start in continuous-time, which we do here. Our main results,
Theorems~\ref{th-mr-ext2} and~\ref{representation},
show that the solutions of~\eqref{spectral}
that are extremal in the max-plus sense can be identified with limits
of ``almost-geodesics'', and that any solution is a max-plus linear combination
of extremal solutions. Our study relies on a compactification of the
state space $\states$, which is similar to the compactification of
metric spaces by horofunctions
(generalised Busemann functions)~\cite{gromov78,rieffel_group}
and to the Martin compactification in potential theory~\cite{martin,dynkin}.

The interest of the notion of almost--geodesic is that it leads to the most
general results since it allows us to avoid making hypotheses on
the Lagrangian or the space $X$ that imply the existence of geodesics.

In Section~\ref{sec-quadratic}, we illustrate these results by
discussing a simple optimal control problem.
We take the state space $\states$ to be $\R^n$, the set
of controls $U$ to be $\R^n$, and $f(x,u):=u$.
The Lagrangian is 
\begin{align}
L(x,u) := -|x|^2 -|u|^2
\enspace, \label{e-lagquad}
\end{align}
where $|\cdot|$ denotes the Euclidean norm.

This optimal control problem
is linear quadratic. In the case where $\lambda=0$,
it is well known that there are
quadratic solutions of~\eqref{spectral}, which can be computed by solving
the algebraic Riccati equation.
These solutions are $\phi(x)=\pm |x|^2$. 
The solution in which one is traditionally interested 
is the ``stable'' one $\phi(x)=-|x|^2$. It yields the feedback
control $\bfu(t)  = \nabla \phi (\bfx(t)) = -2\bfx(t)$,
which pushes the state towards the origin. 
The other quadratic solution yields the opposite control,
sending the state off to infinity, radially. We show however
that there are other solutions, which are non-quadratic.
Amongst them, the extremal solutions, like the unstable
quadratic one, send the state off to infinity, each in
a different direction. 

When $\lambda>0$, there is nothing like the stable quadratic
solution, but the non-quadratic solutions persist.

\section{The max-plus Martin space of a max-plus linear semigroup}
\label{sec-prelim}
Recall that the max-plus semiring $\rmax$ is the set $\R\cup\{-\infty\}$,
equipped with the addition operation $(a,b)\mapsto \max(a,b)$
and multiplication operation $(a,b)\mapsto a+b$.
The completed max-plus semiring $\rmaxb$ is obtained
by adjoining to $\rmax$ a $+\infty$ element, with the convention
that $-\infty$ is absorbing for $(a,b)\mapsto a+b$. 
Although in optimal control applications, the state space $\states$ is typically
a subset of $\R^n$, it is convenient to assume only that $\states$
is an arbitrary set.
To any map $B: \states\times \states \to \rmaxb,\; (x,y)\mapsto B(x,y)$
is associated a max-plus linear self-map from $\rmaxb^\states$ to itself,
$g\mapsto Bg$, defined by:
\begin{align*}
Bg(x) =\sup_{y\in \states} \left(B(x,y)+ g(y)\right) \enspace .
\end{align*}
We say that $(x,y) \mapsto B\of x,y$ is the \new{kernel}
of the map $g\mapsto Bg$. It is uniquely defined.
It is known that a large class of max-plus linear
maps can be represented in this form, see~\cite{DENSITE}.
 
\begin{definition}
We say that a family of operators $(A^t)_{t\geq 0}$ from
$\rmaxb^\states$ to itself is a \new{max-plus linear semigroup with kernel}
if every $A^t$ can be written as
\[
A^tg(x) =\sup_{y\in \states} \left(A^t\of x,y+g(y)\right) \enspace,
\]
for some function $(x,y)\mapsto A^t\of x,y\in \rmaxb$, and if 
\begin{align*}
A^0=I\quad \mrm{and}\quad A^{s+t} = A^s A^t  \enspace,\quad\mrm{ for all } s,t\geq 0 \enspace ,
\end{align*}
where concatenation denotes the composition of operators,
and $I$ denotes the identity operator.
\end{definition}
The semigroup property,
$A^{s+t} = A^s A^t$,
is equivalent to
\begin{align}\label{e-kol}
A^{s+t}\of x,y = \sup_{z\in \states} \left(A^s\of x,z+A^t\of z,y \right)
\enspace .
\end{align}
\begin{example}
\label{exa:laxoleinik}
The Lax-Oleinik semigroup is an example of a max-plus linear
semigroup with kernel: $S^t\of x,y$ is the supremum
of~\eqref{e-reward} over all trajectories $\bfx$ and $\bfu$
satisfying~\eqref{e-dyn}--\eqref{e-init} and
$\bfx(t) = y$. In this case the semigroup property follows from the dynamic programming
principle.
\end{example}
If $(A^t)_{t\geq 0}$ is a max-plus linear semigroup with kernel, then
for all $\lambda\in \R$, 
\[
A^t_\lambda: \qquad g\mapsto -\lambda t + A^tg 
\]
also defines a max-plus linear semigroup, with kernel
$A_{\lambda}^t\of x,y=-\lambda t + A^t\of x,y$. Therefore, when looking for the solutions $\phi$
of~\eqref{spectral}, it is legitimate to assume that $\lambda=0$.
The analogy with potential theory justifies
the following definition.
\begin{definition}
We say that a function $\phi: \states\to \rmax$ is \new{harmonic}
with respect to the semigroup $(A^t)_{t\geq 0}$ if 
\[
A^t\phi =\phi \qquad \mrm{ for all } t \geq 0 \enspace .
\]
We say that $\phi$ is \new{super-harmonic} if 
\[
A^t\phi \leq \phi \qquad \mrm{ for all } t \geq 0 \enspace .
\]
We denote by $\sH$ and $\sS$ the sets of harmonic and super-harmonic functions,
respectively.
\end{definition}
Observe that $\phi$ is required not to take the $+\infty$ value.
We define
\[
A^*:=\sup_{t\geq 0} A^t \enspace.
\]
(The sup is for the pointwise ordering of operators.)
To study harmonic functions, we need to exclude degenerate cases. 
\begin{assumption}\label{irreducible}
We assume that $A^*\of x,y$ is finite for all $x$ and $y$ in $\states$.
\end{assumption}
Observe that $A^*A^*=A^*$, and so
\begin{align}
A^*\of x,z+A^*\of z,y\leq A^*\of x,y 
\; \mrm{ for all } x,y,z\in \states.\label{e-triangular}
\end{align}
This is an analogue of the triangular inequality.
By taking $x=y=z$, we deduce that $A^*\of x,x\leq 0$.
Since $A^*\geq I$, we have $A^*\of x,x\geq 0$, and so
$A^*\of x,x= 0$, for all $x\in \states$.
We fix an arbitrary point $b\in X$, which we call the \new{basepoint}.
\begin{definition}\label{def-martin}
The (max-plus) \new{Martin kernel} of the semigroup $(A^t)_{t\geq 0}$
with respect to the basepoint $b$ is defined by:
\[
K\of x,y = A^*\of x,y- A^*\of b,y \enspace .
\]
The (max-plus) \new{Martin space} $\sM$ of the semigroup $(A^t)_{t\geq 0}$
is the closure in the topology of pointwise convergence
of the set $\sK:=\set{K\of\cdot,y}{y\in \states}\subset \R^\states$.
The (max-plus) \new{Martin boundary} is $\sB=\sM\setminus \sK$.
\end{definition}
{From}~\eqref{e-triangular}, we deduce that $A^*\of b,x+A^*\of x,y\leq A^*\of b,y$,
hence
\[
K\of x,y \leq -A^*\of b,x \enspace ,\qquad \mrm{ for all } x,y\in \states \enspace .
\] 
Similarly,
\[
K\of x,y \geq A^*\of x,b \enspace ,\qquad \mrm{ for all } x,y\in \states \enspace .
\] 
Thus $\sK\subset \prod_{x\in \states} [A^*\of x,b,-A^*\of b,x]$,
and the latter space is compact, by Tychonoff's theorem.
Hence, the closure $\sM$ of $\sK$ is compact.

\begin{remark}
In many applications, the state space $X$ is naturally equipped
with a metric $d$. For instance, when $X\subset \R^n$,
we can take for $d$ the metric induced by any norm on $\R^n$.
Then, the Martin compactification could be defined alternatively
with respect to the topology of uniform convergence
on compact sets, or on bounded sets.
This often leads to the same Martin compactification.
Indeed, the triangular inequality yields
\begin{align}
A^*\of x,y\leq K\of x,z - K\of y,z \leq -A^*\of y,x \enspace .
\label{e-equi}
\end{align}
It is often the case that $|A^*\of x,y | \leq C d(x,y)$ 
for some constant $C>0$. Then,~\eqref{e-equi}
shows that the set $\sK$ is equicontinuous,
and so, by Ascoli's theorem, its closure
in the topology of pointwise convergence
is the same as its closure in the topology of uniform convergence
on compact sets. 
\end{remark}
Any element of $\sK$ is super-harmonic, because $A^tA^*\leq A^*$.
Since a pointwise limit of super-harmonic functions is super-harmonic,
we deduce the following.
\begin{proposition}
Each element of $\sM$ is super-harmonic.
\end{proposition}
As in the discrete time case~\cite{AGW-m}, we have the following.
\begin{proposition}\label{prop-rec}
We have $K\of\cdot,x=K\of\cdot,y$ if and only if
\begin{align}
A^*\of x,y+A^*\of y,x = 0 \enspace .\label{e-rec}
\end{align}
\end{proposition}
In the special case of the Lax-Oleinik semigroup,~\eqref{e-rec} means
that we can find trajectories starting from $x$, passing to $y$,
and returning to $x$, with a reward arbitrarily close to $0$.

Define the equivalence relation $\sim$ on $X$ such that $x\sim y$
if~\eqref{e-rec} holds. Proposition~\ref{prop-rec}
allows us to identify $\sK$ with the quotient $\recclasses:=\states/\sim$.
In particular, when the equivalence classes
of $\sim$ are singletons, the map $x\mapsto K\of \cdot,x$
may be thought of as an ``embedding'' of $X$ in $\sM$,
and the elements of the Martin boundary $\sB=\sM\setminus \sK$
may be thought of as ``boundary points'' of $X$.

For all functions $\xi: \states \to \rmax$ and for all $\eta\in \sM$,
we set:
\begin{align*}
\mu_\xi(\eta )&:= \limsup_{K\of\cdot,x\to \eta} \left( A^*\of b,x+\xi(x) \right)\\
&:= \inf_{W\ni \eta} \sup_{K\of\cdot,x\in W} \left( A^*\of b,x+\xi(x) \right)
\enspace ,
\end{align*}
where the infimum is taken over all neighbourhoods $W$ of $\eta$ in $\sM$.
This map will play the role of the spectral measure in
the Martin representation theorem.

As in the discrete time case, we have the following.
\begin{proposition}\label{u-usc}
Let $u$ be a super-harmonic function.  Then,
$\mu_u(K(\cdot,x)) = A^*\of{b},{x}+ u(x)$ for each $x\in \states$.
\end{proposition}
This result means that the map $\states\to\rmax$, $x\mapsto 
A^*\of{b},{x}+ u(x)$ induces a map on $\recclasses$, or equivalently on
$\sK$, which is upper-semicontinuous in the topology on $\sK$ and
that $\mu_u$ extends this map.

It is convenient to renormalize the operator $A^*$ by 
means of the following transformation. We define,
for all $x,y\in \states$,
\[
A^\natural\of x,y = A^*\of b,x + A^*\of x,y - A^*\of b,y \enspace .
\]
The triangular inequality implies that 
\[
A^\natural\of x,y \leq 0 ,\qquad \text{for all $x,y\in \states$} \enspace .
\]

Proposition~\ref{u-usc} allows us to extend the kernel $A^\natural$
to a kernel $H$ on $\sM\times \sM$ as follows:
\[ H\of \eta,\xi := \mu_\xi(\eta) \enspace .\]
Indeed, 
\[ 
H\of{K\of\cdot,x},{K\of\cdot,y}= A^\natural \of x,y
\enspace .
\]

\begin{definition}
We define the \new{minimal boundary} $\sM^m$ to be the
set of points $\xi$ of $\sM$ that
are (max-plus) harmonic and satisfy $H(\xi,\xi)=0$.
\end{definition}

\section{Main results}

In this section, we state our main results. Many 
of the proofs are similar to those of~\cite{AGW-m}, where
discrete time semigroups are considered. However,
the notion of almost-geodesic appearing there
must be adapted to continuous time.

We use the term path for an arbitrary
map $\gamma$ from a closed interval
of $\R_+$ to $\states$.  
We set, for all $t_0<t_1<\cdots<t_n$ in the interval
of definition of $\gamma$,
\begin{align*}
I(t_0,\dots,t_n;\geo)
   := \sum_{i=0}^{n-1} A^{t_{i+1}-t_{i}}\of{\geo(t_i)},{\geo(t_{i+1})}.
\end{align*}
\begin{definition}
An \emph{almost--geodesic} with parameter $\epsilon$ is a path
$\geo:\rplus\to\states$ such that,
for all $n\in\N$ and $0\leq t_1<t_2<\cdots<t_n$, we have
\begin{align*}
I(t_1,\dots,t_n;\geo) \ge A^*\of\geo(t_1),{\geo(t_n)} - \epsilon.
\end{align*}
\end{definition}

\begin{definition}
A path $\gamma:\rplus\to\states$ is said to be \emph{$\epsilon$-almost--optimal}
with respect to a superharmonic function $h$ if
for all $n\in\N$ and $t_0=0<t_1<t_2<\cdots<t_n$, we have
\begin{align*}
h(\geo(0)) \le \epsilon +
I(t_0,\dots,t_n;\geo) + h(\geo(t_n)).
\end{align*}
\end{definition}

\begin{lemma}\label{optgeo}
Let $\gamma:\rplus\to\states$ be an almost--optimal path with respect
to a superharmonic function.
Then $\gamma$ is an almost--geodesic.
\end{lemma}

\begin{definition}
The reward of a path $\gamma:[s_0,s_1]\to\states$ is
\begin{align*}
I(\geo):=
   \inf_{t_0,\dots,t_n}
   I(t_0,\dots,t_n;\geo),
\end{align*}
where the infimum is taken over all finite increasing sequences
$(t_i)_{i\in\{0,\dots,n\}}$ in $[s_0,s_1]$ with $t_0=s_0$ and $t_n=s_1$.
\end{definition}

If $\geo$ is a path $[s_0,s_1]\to \states$
and $\geo'$ is a path $[s_1,s_2]\to \states$ such that $\geo(s_1)=\geo'(s_1)$, then
we define their concatenation $\gamma\gamma'$ to be the path $[s_0,s_2]\to \states$
that coincides with $\gamma$ on $[s_0,s_1]$ and with $\gamma'$ on $[s_1,s_2]$.

The following result is a consequence of~\eqref{e-kol}.
\begin{lemma}
For all paths $\geo: [s_0,s_1]\to \states$ and $\geo': [s_1,s_2]\to \states$ such that
$\geo(s_1)=\geo'(s_1)$, we have
\[
I(\gamma \gamma') = I(\gamma) + I(\gamma') \enspace .
\]
\end{lemma}

\begin{assumption}\label{assumption}
We assume that for all $t\in\rplus$ and $x$ and $y$ in $\states$,
$A^t\of{x},{y} = \sup_\geo\{I(\geo)\}$, where the supremum is taken over all
paths $\geo:[0,t]\to\states$ from $x$ to $y$.
\end{assumption}
This assumption can be shown to be satisfied in the case when $A^t=S^t$,
the Lax--Oleinik semigroup of example~\ref{exa:laxoleinik}.
However, $I(\gamma)$ will not necessarily be the integral of the Lagrangian
along $\gamma$.

For any path $\alpha:\rplus\to\states$, define
\begin{align*}
J_\alpha(s,t):=
   A^*\of\alpha(0),{\alpha(s)} + I(\alpha|_{[s,t]}) - A^*\of\alpha(0),{\alpha(t)}\enspace .
\end{align*}
Clearly, $J_\alpha$ is always non-positive and
$J_\alpha(s,t)=J_\alpha(s,u)+J_\alpha(u,t)$ for all $s$, $u$, and $t$ in $\R$
with $s\le u\le t$. Using this observation, we can prove the following lemma.
\begin{lemma}\label{tighten}
Let $\geo$ be an $\epsilon$-almost--geodesic,
for some $\epsilon>0$.
Let $y\in\states$ and $\epsilon'>0$.
Then there exists an $\epsilon'$-almost--geodesic $\newgeo$ starting at $y$
and a constant $\Delta\in\R$ such that $\newgeo(t)=\geo(t+\Delta)$
for $t$ large enough.
\end{lemma}
The following lemma is elementary.
\begin{lemma}\label{u-usc2}
Let $u$ be a super-harmonic function.  Then
\begin{align*}
 u=\sup_{w\in \sK} ( \mu_u(w) + w)
=\sup_{w\in \sM} ( \mu_u(w) + w )\enspace .
\end{align*}
\end{lemma}
\begin{prop}\label{semigeo-conv}
If $\geo:\rplus\to\states$ is an almost-geodesic, then $K\of{\cdot},{\geo(t)}$
converges as $t$ tends to infinity to some $w\in \sMin$.
\end{prop}
\begin{proof}
Let $\epsilon>0$.
By Lemma~\ref{tighten}, we lose no generality by assuming
that $\geo$ is an $\epsilon$-almost--geodesic starting at $b$.
So, for all $s$ and $t$ in $\rplus$ with $t\ge s$,
\begin{align*}
A^*\of{b},{\geo(t)}
   &\leq \epsilon + A^*\of{b},{\geo(s)} + A^*\of{\geo(s)},{\geo(t)}\\
   &\leq \epsilon + A^*\of{b},{\geo(t)}
\enspace.
\end{align*}
Therefore
\begin{align}
\label{limKellk}
-\epsilon \le A^*\of{b},{\geo(s)} + K(\geo(s),\geo(t))
\enspace.
\end{align}
Since $\sM$ is compact, it suffices to check that all
convergent subnets of $K(\cdot,\geo(t))$ have
the same limit and that this limit is in $\sMin$.
Let $(\geo({t_d}))_{d\in D}$
and $(\geo({s_e}))_{e\in E}$ be two subnets of $(\geo(t))_{t\in\rplus}$,
such that the nets
\begin{align*}
(K(\cdot,\geo(t_d)))_{d\in D}
\qquad\text{and}\qquad
(K(\cdot,\geo(s_e)))_{e\in E}
\end{align*}
converge respectively to some $w$ and $w'$ in $\sM$.
Applying~\eqref{limKellk} with $s=s_e$ and $t=t_d$, and
taking the limit with respect to $d$,
we obtain $-\epsilon\leq A^*\of{b},{\geo(s_e)}+ w(\geo(s_e))$. Taking now
the limit with respect to $e$,  we get $-\epsilon\leq H(w',w)$.
Since $\epsilon$ is arbitrary and $H(w',w)$ is non-positive,
we obtain $H(w',w)=0$.
{From} Lemma~\ref{u-usc2}, we deduce that $w\geq \mu_{w}(w')+w' =H(w',w)+w'=w'$.
By symmetry, we conclude that $w=w'$, and so $H(w,w)=0$.

Let $y\in\states$ and $\epsilon>0$. By Lemma~\ref{tighten}, there exists
an $\epsilon$-almost--geodesic $\newgeo$ starting at $y$
and such that $K\of{\cdot},{\newgeo(t)}$ converges to $w$.
For all $s$ and $t$ in $\rplus$ with $t\geq s$,
\begin{align*}
A^{s}\of{y},{\newgeo(s)}+A^*\of{\newgeo(s)},{\newgeo(t)}
   \geq  A^*\of{y},{\newgeo(t)} - \epsilon.
\end{align*}
Subtracting $A^*\of{b},{\newgeo(t)}$, 
and letting $t$ tend to infinity, we get
\begin{align*}
A^{s}\of{y},{\newgeo(s)}+w(\newgeo(s)) \ge w(y) - \epsilon.
\end{align*}
Since $\epsilon$ is arbitrary, we have
$w(y)\le \sup_{z\in\states}(A^s\of y, z +w(z))=(A^s w)(y)$.
Since any element of $\sM$ is superharmonic, we get the equality.
\end{proof}

\begin{lemma}
\label{downhill}
Let $h:\states\to\rmax$ be harmonic and let $x\in \states$.
Then for each $\epsilon>0$, there exists a path $\gamma$ starting
at $x$ that is $\epsilon$-almost--optimal with respect to $h$.
\end{lemma}
\begin{proof}
Fix a sequence $\epsilon_n$ in $(0,\infty)$ such that
$\sum_{i=0}^\infty \epsilon_n < \epsilon/2$.
Since $h$ is harmonic, we can construct a sequence $(x_n)_{n\in\N}$ in
$\states$ starting at $x_0:=x$ such that
\begin{align*}
h(x_n) \le \epsilon_n + A^1\of{x_n},{x_{n+1}} + h(x_{n+1})
\qquad\text{for all $n\in\N$}.
\end{align*}
Adding these inequalities, we get
\begin{align}
\label{star1}
h(x) \le
   \sum_{i=0}^{ n-1} A^1\of{x_i},{x_{i+1}} + h(x_{n}) + \sum_{i=0}^{n-1}
\epsilon_i
\quad\text{for all $n\in\N$}.
\end{align}
By assumption~\ref{assumption}, for each $i\in\N$, we can find a path
$\geo_i:[i,i+1]\to\states$ with $\geo_i(i)=x_i$ and $\geo_i(i+1)=x_{i+1}$
such that
\begin{align}
\label{star2}
I(\geo_i)
\geq  A^1\of{x_i},{x_{i+1}} - \epsilon_i
\end{align}
We define $\geo:\rplus\to\states$ to be the
concatenation of these paths, meaning that
$\geo(t):=\geo_{\lfloor t\rfloor}(t)$ for all $t\in\rplus$.
Here $\lfloor t\rfloor$ denotes the largest integer no greater than $t$.
To show that $\geo$ is an $\epsilon$-almost--optimal path with respect to
$h$, consider any increasing sequence of times $t_0,\dots,t_m$
with $t_0:=0$.
Let $n:=\lceil t_m \rceil$ be the least integer no less than $t_m$.
Adding the inequalities obtained from~(\ref{star2}) by letting $i$ run from
$0$ to $n-1$, and then using~(\ref{star1}), we get
\begin{align}
I(\geo|_{[0,n]})&=I(\geo_0\cdots\geo_{n-1})\nonumber \\
&=I(\geo_0)+\cdots +I(\geo_{n-1}) \nonumber\\
   &\geq \sum_{j=0}^{ n-1} A^1\of{x_j},{x_{j+1}} - \sum_{j=0}^{n-1} \epsilon_j 
\nonumber\\
   &\ge h(x)-h(x_n)- 2\sum_{j=0}^{n-1} \epsilon_j.\label{star3}
\end{align}
Using the fact that
$A^{n-t_m}\of{\geo(t_m)},{x_n} \le h(\geo(t_m))-h(x_n)$, we get
\begin{align*}
I(t_0,\dots,t_m;\geo)&= I(t_0,\dots,t_m,n;\geo)-A^{n-t_m}\of{\geo(t_m)},{x_n} 
\\
&\geq  I(\geo|_{[0,n]})+ h(x_n) - h(\geo(t_m)).
\end{align*}
The $\epsilon$-almost--optimality of $\geo$ follows on combining this
with~(\ref{star3}).
\end{proof}
The following is a max-plus analogue of the Martin representation theorem
in potential theory. A discrete time version was
given in~\cite{AGW-m}.
\begin{theorem}\label{representation}
A function $h:\states\to\rmax$ is harmonic if and only if it can be written
\begin{align}
\label{star4}
h=\sup_{w\in \sMin} (\nu(w) + w) \enspace ,
\end{align}
where $\nu$ is some upper semicontinuous function from $\sMin$ to $\rmax$.
Moreover, $\mu_h$ is the greatest $\nu$ satisfying this equation.
\end{theorem}
\begin{proof}
Let $h$ be harmonic. By Lemma~\ref{u-usc2},
\begin{align*}
h=\sup_{w\in  \sM} (\mu_h(w) + w)
   \geq \sup_{w\in \sMin}( \mu_h(w)+ w) \enspace .
\end{align*}
To show the opposite inequality, fix $x\in \states$ such that $h(x)>-\infty$.
Let $\epsilon>0$. By Lemma~\ref{downhill}, there exists an
$\epsilon$-almost--optimal path $\geo:\rplus\to\states$ with respect to $h$
starting at $x$. By Lemma~\ref{optgeo}, this path is an almost--geodesic
and therefore, by Lemma~\ref{semigeo-conv}, $K\of\cdot,{\geo(t)}$
converges to a point $\xi$ of $\sMin$. Also, the $\epsilon$-almost--optimality
of $\geo$ implies that
\begin{align*}
h(x) &\le \epsilon + A^*\of{x},{\geo(t)} + h(\geo(t)) \\
     &= \epsilon + K\of{x},{\geo(t)} + A^*\of{b},{\geo(t)} + h(\geo(t)) \enspace.
\end{align*}
Letting $t$ tend to infinity and using the fact that $\epsilon$ is arbitrary,
we see that $h(x) \le \xi(x)+\mu_h(\xi)$.
We have thus established that
\begin{align*}
h(x) \le  \sup_{w\in \sMin}( \mu_h(w)+ w) \enspace .
\end{align*}

Now let $\nu:\sMin\to\rmax$ be any function satisfying~(\ref{star4}) and let
$\xi\in\sMin$. Then $h\ge \nu(\xi)+\xi$. So
\begin{align*}
\mu_h(\xi) &= \limsup_{K\of\cdot, x\to \xi}(A^*\of{b},{x}+h(x))\\
    &\ge \limsup_{K\of\cdot, x\to \xi}(A^*\of{b},{x}+\xi(x))+\nu(\xi) \\
    &= H(\xi,\xi)+\nu(\xi) \\
    &= \nu(\xi) \enspace .
\end{align*}
Therefore $\nu\le \mu_h$.
\end{proof}

The remaining results can be proved by adapting
the arguments of the corresponding results in~\cite{AGW-m}.
\begin{lemma}\label{prop-xi}\sloppy
Let $\xi\in \sM$ be such that $H(\xi,\xi)=0$ and suppose that
$\xi$ can be written in the form $\xi=\sup_{w\in \sM}( \nu(w)+ w)$,
where $\nu:\sM\to\rmax$ is upper semicontinuous.
Then, there exists $w\in \sM$ such that $\xi=\nu(w)+ w$.
\end{lemma}

\begin{definition}
Let $V$ be a set of functions from $X$ to $\R$ that is closed under
maximisation and the addition of a constant.
A function $f$ from $X$ to $\R$ is said to be
\emph{normalised} if $f(b)=0$,
and is said to be an \emph{extremal generator} of $V$
if $f=\max(u,v)$ with
$u$ and $v$ in $V$ implies that $f$ is equal to either $u$ or $v$.
\end{definition}
\begin{corollary}\label{mr-ext}
Every element $\xi$ of $\sM$ satisfying $H(\xi,\xi)=0$
is a normalised extremal generator of $\sS$.
\end{corollary}

\begin{theorem}\label{th-mr-ext2}
The normalised extremal generators of $\sH$ are precisely the elements
of $\sMin$.
\end{theorem}

\begin{prop}
Assume that $\sM$ is first countable. For all $\xi\in \minspace$, there
exists an almost--geodesic converging to $\xi$.
\end{prop}

\section{The Horoboundary of a Linear Quadratic Control Problem}
\label{sec-quadratic}
We next study the eigenproblem~\eqref{spectral}, for the Lax-Oleinik
semigroup $(A^t)_{t\geq 0}$ corresponding to the Lagrangian~\eqref{e-lagquad}. 

Observe first that $A^t\phi(0)\geq \phi(0)$, since the zero control makes
the trajectory stay at the origin with a zero cost. 
Hence, any eigenvalue $\lambda$ of the semigroup must be nonnegative.

The normalised semigroup $A^T_\lambda=-\lambda T+ A^T$
corresponds to the reward:
\begin{align*}
-\int_0^T (|\bfx(t)|^2+|\dot \bfx(t)|^2+\lambda)\,{\rm d}t,
\end{align*}
meaning that $A^T_\lambda\of x,y$ is the maximal value of this reward
over all curves $\bfx$ joining $x$ to $y$.

One may regard the optimisation of this reward
as a physics problem in which a particle is subjected
to a force away from the origin proportional in strength to its distance
from the origin.
Indeed, the Euler equation gives that optimal paths satisfy $\Ddot \bfx(t)=\bfx(t)$
for all $t\in\R$.
The general solution to this differential equation is $\bfx(t)= W e^t+Z e^{-t}$,
where $W$ and $Z$ are constant vectors.
If the position is $x$ at time $0$ and $y$ at time $T$, then
\begin{align*}
W=\frac{y-e^{-T}x}{e^{T}-e^{-T}}
\qquad\text{and}\qquad\
Z=\frac{e^{T}x-y}{e^{T}-e^{-T}}.
\end{align*}
Doing the integration, we find that the action of the optimal path is
\begin{align}
\label{pathaction}
-\frac{(|x|^2+|y|^2)\cosh T -2x\cdot y}{\sinh T}-\lambda T \enspace .
\end{align}

Now we must maximise this expression over $T\ge 0$.
We will first do the case when $\lambda=0$ since the formulae are more
manageable. If $x\cdot y \le 0$, then the supremum occurs at $T=+\infty$,
where the function above tends to the limit $A^*\of{x},{y}=-|x|^2-|y|^2$.
On the other hand, if $x\cdot y > 0$, then the supremum is attained
at finite $T$. Differentiating the expression above and setting it equal to
zero, we obtain a linear equation in $\cosh T$, with solution
\begin{align*}
\cosh T = \frac{|x|^2+|y|^2}{2x\cdot y}.
\end{align*}
{From} this, we can easily calculate $\sinh T$ as well.
Substituting these values into~(\ref{pathaction}) and simplifying,
we obtain
\begin{align*}
A^*\of{x},{y}=-|x-y||x+y|
\qquad\text{if $x\cdot y>0$}.
\end{align*}
\addtolength{\textheight}{-0.9cm}
To calculate the boundary points we take the limit of $A^*\of{x},{y}-A^*\of{0},{y}$
as $|y|$ tends to infinity, keeping $n:= y/|y|$ constant.
One expands the expression in powers of $|y|$ and discards those that tend
to zero. The resulting boundary point is the function:
\begin{align*}
h_n(x) = \begin{cases}
            -|x|^2+2(x\cdot n)^2, & \text{if $x\cdot n>0$}, \\
            -|x|^2, & \text{otherwise}.
   \end{cases}
\end{align*}

\begin{figure}
\begin{center}
\includegraphics{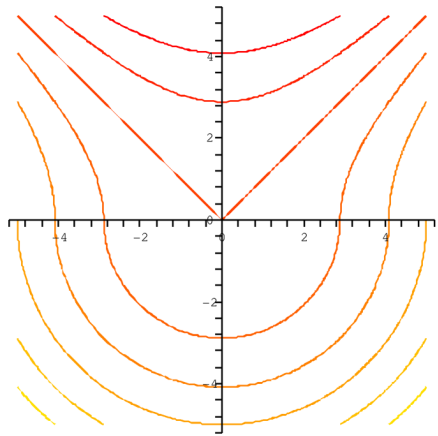}\end{center}
\caption{Horospheres of the LQ model when $\lambda=0$.}
\label{fig1}
\end{figure}

Figure~\ref{fig1} shows the ``horospheres'' of $h_n$
in dimension two, when $n=(0,1)$. These are the sets
of the form $\set{x\in \R^2}{h_n(x)=\mrm{const}.}$.

Now we turn to the case when $\lambda>0$.
To find the maximum over $T\ge 0$ of~(\ref{pathaction}), we again
differentiate it and set it equal to zero. This time we obtain a quadratic
equation in $\cosh T$, with solutions
\begin{align*}
\lambda\cosh T
   = -x\cdot y \pm \sqrt{(x\cdot y)^2 + \lambda^2 + \lambda(|x|^2+|y|^2)}.
\end{align*}
It is of course the positive solution we want.
We next calculate $\sinh T$ and $T$.
Substituting them into~(\ref{pathaction}), we obtain an expression
for $A^*\of{x},{y}$, which, unfortunately, is now rather complicated.
However the special case when $x=0$ is not too bad:
\begin{align*}
A^*\of{0},{y} :=
-  |y|\sqrt{\lambda\!+\!y^2}-\lambda\log(\sqrt{\lambda\!+\!y^2}+|y|)+\lambda\log\lambda.
\end{align*}
To calculate the boundary points we again take the limit of
$A^*\of{x},{y}-A^*\of{0},{y}$ as $|y|$ tends to infinity, keeping $n:= y/|y|$ constant.
The best way of calculating this limit is to treat separately the terms
involving a logarithm from the rest. Again one should expand the expression
as powers of $|y|$, discarding those that tend to zero.
After much simplification, one arrives at the result
\begin{align*}
h_n(x) =
-   \lambda\frac{|x|^2}{R^2}
      + x\cdot n\frac{\lambda+ 2|x|^2}{R}
      - \lambda\log\frac{R}{\sqrt\lambda},
\end{align*}
where $R:=\sqrt{(x\cdot n)^2+ \lambda}-x\cdot n$.

\begin{figure}
\begin{center}\includegraphics{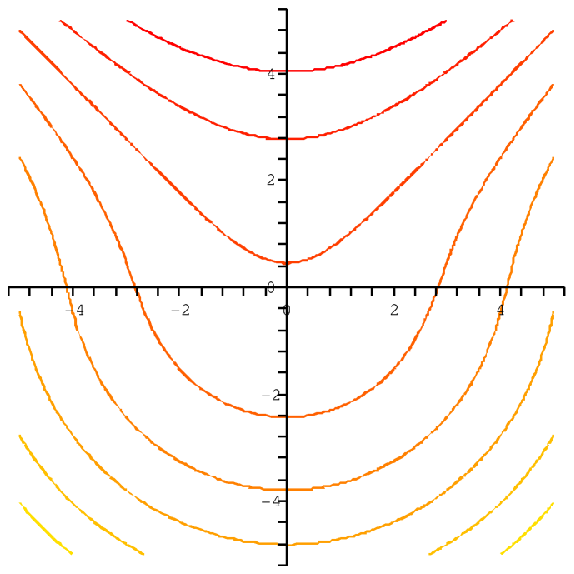}\end{center}
\caption{Horospheres of the LQ model when $\lambda=1$.}
\label{fig2}
\end{figure}

Figure~\ref{fig2} shows the horospheres of $h_n$
when $n=(0,1)$, again in dimension two.

\bibliographystyle{plain}
\bibliography{cdc06b}
\end{document}